\documentclass[12pt]{amsart}
\usepackage{amsthm}
\usepackage{comment}
\usepackage{enumerate}
\usepackage{amsmath}
\usepackage{float}
\usepackage{amsfonts}
\newtheorem{theorem}{Theorem}[section]
\newtheorem{lemma}[theorem]{Lemma}

\newtheorem{corollary}[theorem]{Corollary}
\newtheorem{notation}[theorem]{Notation}
\newtheorem{conjecture}[theorem]{Conjecture}

\newtheorem{definition}[theorem]{Definition}

\newtheorem{remark}[theorem]{Remark}

\restylefloat{figure}
\restylefloat{table}

\def\R{{\rm Res}}

\def\OO{{\mathcal O}}

\def\CC{{\mathbb C}}

\def\NN{{\mathbb N}}

\newcommand{\Cline}[1]{\vspace{-21 pt} \\ \cline{#1} \\ \vspace{-21 pt} \\}

\begin{document}
\title{On the height of the Sylvester Resultant}
\author{Carlos D'Andrea}
\thanks{Research of C. D'Andrea supported by 
 the Miller Institute for Basic Research in Science, University of California, Berkeley}
\address{Department of Mathematics, University of California, Berkeley,
    USA, 94720-3840, 970 Evans Hall}
\email{cdandrea@math.berkeley.edu}

\author{Kevin G. Hare}
\thanks{Research of K.G. Hare supported, in part by NSERC of Canada}
\address{Department of Pure Mathematics, University of Waterloo, Waterloo,
    Ontario, Canada,  N2L 3G1}
\email{kghare@cecm.sfu.ca}

\begin{abstract}
Let $n$ be a positive integer. 
We consider the Sylvester resultant of $f$ and $g,$ where $f$ is a generic polynomial of degree $2$ or $3$ and $g$ is a generic polynomial of degree $n.$
If $f$ is a quadratic polynomial, we find the resultant's height. If
$f$ is a cubic polynomial, we find tight asymptotics for the resultant's
    height.
\end{abstract}

\date{\today}
\maketitle

\section{Introduction}
\label{sec:intro}
Let $m$ and $n$ be positive integers, $f$ and $g$ be generic univariate polynomials of degrees $m$ and $n$ respectively:
\begin{equation}\label{syste}
\begin{array}{ccc}
f(x)&:=&f_0+f_1x+\cdots+f_mx^m, \\
g(x)&:=&g_0+g_1x+\cdots+g_nx^n.
\end{array}
\end{equation}
Here, $f_i,\,g_j$ are new variables. The Sylvester resultant of $f$ and $g$ is the determinant of the following square matrix of order $m+n:$
\begin{equation}\label{def}
\R(f,g) := \det\left[\begin{array}{ccccccc}
f_0    &        &        &        & g_0    &        &          \\
f_1    & f_0    &        &        & g_1    & \ddots &          \\
\vdots & \vdots & \ddots &        & \vdots & \ddots & g_0      \\
f_m    & f_{m-1}&        & f_0    & g_{n-1}&        & g_1      \\
       & f_m    & \ddots & \vdots & g_n    & \ddots & \vdots   \\
       &        & \ddots & f_{m-1}&        & \ddots & g_{n-1}  \\
       &        &        & f_m    &        &        & g_n      \\
\end{array}
\right],
\end{equation}
where the first $n$ columns contain coefficients of $f$ and the last $m$ contain coefficients of $g.$ 

From the definition, it is very easy to see that $\R(f,g)$ is a homogeneous polynomial in the variables $f_i$ and $g_j.$ 
Further $\R(f,g)$ is homogeneous in
each group of variables, having degree $n$ in the $f_i$'s, and $m$ in the $g_j$'s. It is not hard to see that the resultant is 
$\omega$-homogeneous of ``degree'' $mn,$ where $\omega=(0,1,\cdots,n,0,1,\cdots,m).$ 
This means that if a monomial $f_0^{\alpha_0}\cdots f_m^{\alpha_m} g_0^{\beta_0}\cdots g_n^{\beta_n}$ 
appears with nonzero coefficient in the expansion of $\R(f,g)$ then $\sum_{i=1}^mi\alpha_i+\sum_{j=1}^nj\beta_j=mn$ (see \cite[Theorem $6.1$]{stu}).

Resultants are widely used as a tool for polynomial equation solving, this has sparked a lot interest in their computation (see e.g. \cite{clo,clo2,GKZ}).
The absolute height of a polynomial $g=\sum_{\alpha}c_\alpha U^\alpha\in\CC[U_1,\cdots,U_p]$ is defined as
$H(g):=\max\{|c_\alpha|,\,\alpha\in\NN^p\}.$ In this paper we will be concerned with the computation of the height of $\R(f,g).$

The sharpest upper bound for the height was given in \cite[Theorem $1.1$]{som}, where it is shown that 
$H\left(\R(f,g)\right)\leq (m+1)^n\,(n+1)^m.$ 
Previous upper bounds were given in 
\cite{BGS,KPS,phi,roj,som1}, for more general resultants which include
$R(f,g).$

However, up to now there have been no known exact expressions for $H(\R(f,g))$, for any non-trivial cases. 
We only know the exact value of the coefficients of the resultant for extremal monomials with respect to a generic weight, and they are equal to $\pm1$ 
(see \cite[Corollary $3.1$]{stu}).

The purpose of this paper is to give non-trivial estimates on the height of the resultant for polynomials $f$ of low degree.

\subsection{Quadratic polynomials} 
In the case $m=2$, we get an exact solution for the height of $\R(f,g)$
    in terms of an integer number $A_n$.
To define $A_n$, first consider $p_n(z):=(n-2z+1)(n-2z+2)-z(n-z).$
It is easy to see that if $n \geq 3$, then $p_n(0)>0$ and 
    $p_n\left(\frac{n}{2}\right)<0.$
As $p_n(z)$ is a quadratic polynomial in $z,$ we define, for 
$n\geq3,$  $r_n$ as the unique root of $p_n(z)$ lying in 
    $\left[0,\frac{n}2\right].$ 
Set $A_n:=\lfloor r_n\rfloor$, the floor of $r_n$.
In Table \ref{tab:An} page \pageref{tab:An}, we have listed some values of 
    $A_n.$

\begin{theorem}\label{thm:mt1}
Let $n\geq3.$ The coefficient of highest absolute value in the expansion of $\R(f_0+f_1x+f_2x^2,g)$ 
is the coefficient corresponding to $g_0g_nf_0^{A_n}f_1^{n-2A_n}f_2^{A_n}.$
Moreover, 
\begin{eqnarray*}
H\left(\R(f_0+f_1x+f_2x^2,g)\right) 
& = & H\left(\R(f_0+f_1x+f_2x^2,g_0 + g_n x^n)\right) \\
& = & n\frac{(n-A_n-1)!}{(n-2A_n)!A_n!}.
\end{eqnarray*}
\end{theorem}

\begin{remark}
As $A_n<\frac{n}{2},$ it turns out that $(n-2A_n)\geq0.$
\end{remark}

Before we give the next result, we must introduce some notation.

\begin{notation}
Let $\alpha(n)$ be a positive sequence.
We say that a sequence $\beta(n) = \OO(\alpha(n))$ if there exist positive constants
    $c_1, c_2$ and $N$ such that for all $n > N$ we have
    $c_1 \alpha(n) \leq \beta(n) \leq c_2 \alpha(n)$
\end{notation}

Based on Theorem \ref{thm:mt1} we get
\begin{corollary}
\label{cor:alpha2}
Let $\alpha \approx 1.6180$ be the positive root of 
    $x^2-x-1$ and 
    $\beta \approx 2.3644$ be the positive root of $4 x^4 - 125.$
Then 
   \[H\left(\R(f_0+f_1x+f_2x^2,g)\right) = \frac{\beta}{\sqrt{n \pi}}\alpha^n
    - \OO\left(\frac{\alpha^n}{n \sqrt{n}}\right)\]
\end{corollary}

\subsection{Cubic polynomials}

In the case  $m = 3$, we get a tight bound for the height.
In particular, we get the following:
\begin{theorem}\label{thm:mt2}
Let $\beta \approx 8.13488$ be the real root of $x^3-18 x^2+110 x-242$, and
    $\alpha \approx 1.83928$ be the real root of $x^3-x^2-x-1$.
Let $H(n) := H\left(\R(f_0+f_1x+f_2x^2+f_3x^3,g)\right)$ where
    $g$ is of degree $n$.
Then
\begin{equation}\label{cucu}
H(n) = \frac{\beta}{\pi n}\alpha^n - \OO\left(\frac{\alpha^n}{n^2}\right)
\end{equation}
\end{theorem}

\subsection{Organization of paper}

Section \ref{sec:bi} gives a proof of Theorem \ref{thm:mt1}
    and Corollary \ref{cor:alpha2}.
A proof of Theorem \ref{thm:mt2} is given in Section \ref{sec:tri}.
Section \ref{sec:conc} gives some conclusions, conjectures and list
    some open questions.
Finally, Section \ref{sec:append} contains a number of different
    tables which are referred to throughout this paper.

\section{Quadratic polynomials}
\label{sec:bi}

\begin{proof}[Proof of Theorem \ref{thm:mt1}]
The proof will be made by induction on $n.$ 
For this section, denote with $H(n)$ the height of the resultant of 
a degree-two generic polynomial $f$ and a generic polynomial $g$ of degree $n.$

For $n=3,$ an explicit computation shows that
\begin{itemize}
\item $A_3=1,$
\item $H(3) = 3,$ and this is the coefficient of $g_0g_3f_0f_1f_2.$
\end{itemize}

Suppose now $n>3.$ As the degree of ${\R}(f,g)$ in the $g_j$'s is $2,$ we will first consider two special cases:
\begin{itemize}
\item If we pick a term in the expansion of ${\R(f,g)}$ which is not a multiple of $g_0,$ this term will appear in the expansion of
$${\R}(f,g_nx^n+\dots+g_1x)=\pm f_0{\R}(f,g_nx^{n-1}+\dots+g_1),$$
and by the inductive hypothesis, all the coefficients of this expansion are bounded by $H(n-1).$
\item If we pick a term in the expansion of ${\R(f,g)}$ which is not a multiple of $g_n,$ this term will appear in the expansion of 
$${\R}(f,g)=\pm f_2 {\R}(f, g_{n-1}x^{n-1}+\dots+g_0),$$
and reasoning as in the previous case, all the coefficients in this case will be bounded by $H(n-1).$
\end{itemize}
In order to conclude, we have to bound all the coefficients which appear in ${\R(f,g)}$ which are coefficients to a monomial of the form $g_0g_nf_0^a f_1^b f_2^c,$ for some $a, b$ and $c$, and compare this bound with $H(n-1).$

Without loss of generality we compute 
    ${\R} (f_2x^2+f_1 x+f_0, g_nx^n+g_0).$ 
Moreover, we can also set $g_n:=f_2:=1.$ 
Let $f(x)=(x-x_1)(x-x_2).$ Then,
\begin{equation}
\begin{array}{rcl}
{\R}(f,g)& =& \pm({x_1}^n+g_0)({x_2}^n+g_0) \\
         & =& \pm\left({(x_1x_2)}^n+({x_1}^n+{x_2}^n)g_0+g_0\right).
\end{array}
\label{uno}
\end{equation}
In order to write the right-hand side of (\ref{uno}) in terms of $f_1,f_0,$ we apply the classical Girard formulas (see for instance \cite[Chapter $4$ F]{GKZ}):
\begin{equation}\label{dos}
{x_1}^n+{x_2}^n = (-1)^nn\sum_{i_1+2i_0=n}{(-1)^{2i_1+i_0}\frac{(i_1+i_0-1)!}{i_1!i_0!}{f_1}^{i_1}{f_0}^{i_0}}.
\end{equation}
So, we have to maximize $\frac{(i_1+i_0-1)!}{i_1!i_0!}$ subject to 
the condition $i_1+2i_0=n.$
Set $z:=i_0,$ then $i_1=n-2z,$ and we have study the behaviour of the function
$$P(z):=\frac{(n-z-1)!}{(n-2z)!z!},\ \mbox{for}\,z=0,1,\dots,
    \left\lfloor\frac{n}{2}\right\rfloor.$$

As $$P(z)-P(z-1)=\frac{(n-z-1)!}{(n-2z+2)!z!} p_n(z),$$
and due to the fact that $p_n(z)$ is a quadratic equation having $r_n$ as the unique root in the interval $[0,\frac{n}{2}],$ we have
\begin{itemize}
\item $P$ is increasing for $z=0,1,\dots,A_n.$
\item $P$ decreases for $z=A_n, A_n+1,\dots,\left\lfloor\frac{n}2\right\rfloor$.
\end{itemize}
Hence, the maximum of $P$ is attained when $z=A_n,$ and $H(n) = n P(A_n)$ 
    because of (\ref{uno}) and (\ref{dos}).

In order to conclude, we only have to prove that $H(n)>H(n-1).$
As $H(n-1)=(n-1)\frac{(n-A_{n-1}-2)!}{(n-1-2A_{n-1})!A_{n-1}!},$
and 
\begin{equation}\label{tres}
H(n)\geq n\frac{(n-A_{n-1}-1)!}{(n-2A_{n-1})!A_{n-1}!},
\end{equation}
it is easy to check that the right-hand-side of (\ref{tres}) is bigger than $H(n-1)$ if and only if $n\geq 3.$ 
\end{proof}

From here, we can prove Corollary \ref{cor:alpha2}:

\begin{proof}[Proof of Corollary \ref{cor:alpha2}]
By noticing that $r_n = \frac{6 + 5 n - \sqrt{5 n^2 - 4}}{10}$, 
    we get
\[\lim_{n\rightarrow\infty} \frac{A_n}{n} = \frac{5 - \sqrt{5}}{10}.\]
Thus for large $n$ we get
\begin{eqnarray*}
n \frac{(n - A_n - 1)!}{(n-2 A_n)!A_n!} 
& = & n \frac{\Gamma(n - A_n)}{\Gamma(n - 2 A_n + 1)\Gamma(A_n + 1)} \\
& = & \frac{n \Gamma(n-A_n)}{(n-2 A_n) A_n \Gamma(n-2 A_n) \Gamma(A_n+1)} \\
& = & \frac{n^2}{(n-2 A_n) A_n}\times \frac{\Gamma(n-A_n)}{n \Gamma(n-2 A_n) 
                     \Gamma(A_n+1)} \\
\end{eqnarray*}

From the comment above, we see that the first fraction will approach
    to $\frac{5 (1 + \sqrt{5})}{2}$.
This then gives us

\begin{eqnarray*}
& \approx & \frac{5 (1+\sqrt{5})}{2}\frac{\Gamma(n/2 + n \sqrt{5}/10)}
       {n \Gamma(n \sqrt{5}/5) \Gamma(n/2 - n \sqrt{5}/10) } \\
& = & \frac{\beta}{\sqrt{\pi n}}\alpha^n
       - \OO\left(\frac{\alpha^n}{n^{3/2}}\right) \\
\end{eqnarray*}
which gives the desired result.
The last line of this inequality was derived with the help of Maple.

Here we ignored a number of problems that occur with respect to errors in approximation.
These are done in the same way that they are done for the proof of Theorem \ref{thm:tight}.
\end{proof}

\section{Cubic polynomials}
\label{sec:tri}

By an argument similar to Theorem \ref{thm:mt1}, if $H(n) > H(n-1)$ then both $g_n$ and $g_0$ must divide the terms which gives rise to $H(n)$. 
We will see that this holds for $n\gg0.$
We have then that three $g_i$'s must divide each of the terms of ${\rm Res}(f,g)$ and two of them are known if $H(n)>H(n-1)$ ($g_n$ and $g_0$).
This gives rise to the definitions
\begin{definition}
\label{defn:Hm}
Define $H_l(m,k,k',m')$ to be the coefficient 
    of $f_0^m f_1^k f_2^{k'} f_3^{m'}g_0 g_l g_n$
    in $\R(f, g)$.
\end{definition}

\begin{definition}
\label{defn:H}
Define \[H_l(n) = \max_{m+k+k'+m'=n} \left|H_l(m,k,k',m')\right|.\]
\end{definition}

The main results of the paper will be derived by being able to write
    $H_l(m,k,k',m')$ in terms of some auxiliary functions $F(m,k,k',m')$ which are defined as follows:

\begin{definition}
\label{defn:F}
Define $F(m,k,k',m')$ to be the number of occurrences of 
    $f_0^m f_1^k f_2^{k'} f_3^{m'}$ in the determinant of the matrix
\[\left[\begin{array}{ccccccc}
f_2 & f_1 & f_0 \\
f_3 & f_2 & f_1 & f_0 \\
    & f_3 & f_2 & f_1 & f_0 \\
    &     & \ddots & \ddots & \ddots & \ddots \\
    &     &     & f_3 & f_2 & f_1 & f_0 \\
    &     &     &     & f_3 & f_2 & f_1 \\
    &     &     &     &     & f_3 & f_2 \\
\end{array}\right]\]
of dimension $m+k+k'+m' \geq 1$.
For $m+k+k'+m'=1$ or $2$ the determinant would be of the matrices
    $[f_2]$ and 
    $\left[\begin{array}{cc}f_2 & f_1 \\ f_3 & f_2 \end{array}\right]$
    respectively.

For convenience we define $F(0,0,0,0) = 1$.
\end{definition}

For example for $m+k+k'+m'=3$ we have 
\[ \det\left[
\begin{array}{ccc}
f_2 & f_1 & f_0 \\
f_3 & f_2 & f_1 \\
 0  & f_3 & f_2 \\
\end{array}\right] = f_2^3 - 2 f_1 f_2 f_3 + f_0 f_3^2
\]
Thus we see that $F(1,0,0,2) = 1, F(0,1,1,1) = -2$ and $F(0,0,3,0) = 1$.

\begin{lemma}
\label{lem:rec}
$F(m, k, k', m')$ satisfies the recurrence relation
    \begin{eqnarray*}
F(m, k, k', m') & = & F(m, k, k'-1, m') - F(m, k-1, k', m'-1)  \\ && 
                       + F(m-1, k, k', m'-2) 
\end{eqnarray*}
with $F(0,0,0,0) = 1$ and $F(m,k,k',m') =0$ if any of $m, k, k'$ or $m' < 0$
\end{lemma}

\begin{proof}
The recurrence follows by considering the three possibilities from the 
    first row.
\begin{eqnarray*}
&& \left[\begin{array}{ccccccc}
\fbox{$f_2$} & f_1 & f_0 \\
\Cline{2-7}
f_3 & f_2 & f_1 & f_0 \\
\Cline{1-1}
    & f_3 & f_2 & f_1 & f_0 \\
    &     & \ddots & \ddots & \ddots & \ddots \\
    &     &     & f_3 & f_2 & f_1 & f_0 \\
    &     &     &     & f_3 & f_2 & f_1 \\
    &     &     &     &     & f_3 & f_2 \\
\end{array}\right], 
\left[\begin{array}{ccccccc}
f_2 & \fbox{$f_1$} & f_0 \\
\Cline{1-1} 
\Cline{3-7}
\fbox{$f_3$}& f_2 & f_1 & f_0 \\
\Cline{2-7}
    & f_3 & f_2 & f_1 & f_0 \\
\Cline{1-2}
    &     & \ddots & \ddots & \ddots & \ddots \\
    &     &     & f_3 & f_2 & f_1 & f_0 \\
    &     &     &     & f_3 & f_2 & f_1 \\
    &     &     &     &     & f_3 & f_2 \\
\end{array}\right] \\
 && \left[\begin{array}{cccccccc}
f_2 & f_1 & \fbox{$f_0$} \\
\Cline{1-2}
\Cline{4-8}
\fbox{$f_3$} & f_2 & f_1 & f_0 \\
\Cline{2-8}
    & \fbox{$f_3$} & f_2 & f_1 & f_0 \\
\Cline{1-1} 
\Cline{3-8}
    &     & f_3 & f_2 & f_1 & f_0 \\
\Cline{1-3}
    &     &     & \ddots & \ddots & \ddots & \ddots \\
    &     &     &     & f_3 & f_2 & f_1 & f_0 \\
    &     &     &     &     & f_3 & f_2 & f_1 \\
    &     &     &     &     &     & f_3 & f_2 \\
\end{array}\right].
\end{eqnarray*}
\end{proof}
    
By induction we will prove the following lemma,
whose statement was first discovered experimentally via \cite{SloaneOnline}.
\begin{lemma} 
\label{lem:binomial}
If $m' = 2 m + k$, then:
\begin{equation}
F(m, k, k', k+2 m) = (-1)^k {m+k \choose k} {k' + k + m\choose k+m} 
\label{eq:choose}
\end{equation}
If $m' \neq 2 m + k$ then $F(m, k, k', m') = 0$.
\end{lemma}

\begin{proof}
By examining the recurrence relation, we see that $F(m,k,k',m') = 0$ if 
    $m' \neq 2 m + k$.

Equation (\ref{eq:choose}) is true for $m + k + k' = 1$ by some simple 
    calculations.
So we have that 
\begin{eqnarray*}
& & F(m, k, k', k+2 m)  \\
    & = & F(m, k, k'-1, k + 2m) - F(m, k-1, k', k+2m-1) \\ 
    &   & +F(m-1, k, k', k+2m-2) \\
    & = & (-1)^k{m+k\choose k}{k'-1+k+m\choose k+m}  \\ &&
          - (-1)^{k-1} {m+k-1 \choose k-1}{k' + k-1 + m \choose k+m-1} \\ &&
          + (-1)^k     {m+k-1 \choose k}{k'+k+m-1 \choose k+m-1} \\
    & = &   (-1)^k \left({m+k\choose k}     {k'-1+k+m\choose k+m} \right. \\ & &
          \left.+ {k' + k-1 + m \choose k+m-1} 
            \left({m+k-1 \choose k-1} + {m+k-1 \choose k}\right)\right) \\
    & = &   (-1)^k \left({m+k\choose k}\left({k'-1+k+m\choose k+m}  
          + {k' + k-1 + m \choose k+m-1}\right)\right) \\
    & = &   (-1)^k {m+k\choose k}{k'+k+m\choose k+m}  
\end{eqnarray*}
and the result follows by induction.
\end{proof}

\begin{theorem}
\label{thm:main}
Let $F$ be as in Definition \ref{defn:F}.
Then
\begin{eqnarray*}
H_0(m,k,k',m') 
    & = & F(m-1,k,k',m'-2)- F(m,k,k'-1,m') \\
    &   & +2  F(m,k,k',m') \\
    & = & (-1)^k (3 m + 2k + k')\frac{(m+k+k'-1)!} {k! m! k'!}.
\end{eqnarray*}
\end{theorem}

The value of $H_l(m,k,k',m')$ is given in Table \ref{tab:recurrence}
    (page \pageref{tab:recurrence}) for $l$ from 0 to 5.
We will provide only the proof for $H_0(m,k,k',m')$ here.
The other cases listed in Table \ref{tab:recurrence} are similar.
Code which automates this process is available upon request.

For all $l$, we can also write $H_l(m,k,k',m')$ as a sum of 
    various $F$.
Instead of three cases, we tend to get six, depending
    on which column the $g_0$, the $g_l$ and the $g_n$ are taken from.
In each of these cases we get a finite number of ways to account for the
    terms above the $g_l$ term, and below the $g_n$ term.
The terms between the $g_l$ and the $g_n$ can be accounted for with $F$ 
    functions.
So each of these finite number of ways will account for some 
    $F(m-?, k-?, k'-?, m'-?)$ which will then be taken into the final 
    sum.

\begin{proof}[Proof of Theorem \ref{thm:main}]
The second statement of the theorem follows directly from Lemma 
    \ref{lem:binomial}, so it suffices to prove the first statement.

We notice that there are three different ways in which we can get 
    $g_0 g_0 g_n$ as a factor.
We will do each case separately.
\begin{enumerate}[{Case} 1:]
\item
\[\left[
\begin{array}{cccccccccc}
f_0 &     &     &     &        &      &     & \fbox{$g_0$} &    & \\ 
\Cline{1-7} \Cline{9-10}
f_1 & f_0 &     &     &        &      &     & g_1 & \fbox{$g_0$}& \\
\Cline{1-8} \Cline{10-10}
f_2 & f_1 & f_0 &     &        &      &     & g_2 & g_1 & g_0 \\
f_3 & f_2 & f_1 & \ddots &     &      &     & \vdots & \vdots & \vdots \\
    & f_3 & f_2 & \ddots & f_0 &      &     & \vdots & \vdots & \vdots \\ 
    &     & f_3 & \ddots & f_1 & f_0  &     & \vdots & \vdots & \vdots \\ 
    &     &     & \ddots & f_2 & f_1  &     &  g_n   & g_{n-1}&  g_{n-2} \\ 
    &     &     &        & f_3 & f_2  &     &        &  g_n   &  g_{n-1} \\ 
    &     &     &        &     & f_3  &     &        &        & \fbox{$g_n$} \\
\Cline{1-9}
\end{array}\right]\]

So we get that this case contributes \( F(m, k, k', m').\)

\item
\[\left[
\begin{array}{cccccccccc}
f_0 &     &     &     &        &      &     & \fbox{$g_0$} &    & \\ 
\Cline{1-7} \Cline{9-10}
f_1 & f_0 &     &     &        &      &     & g_1 & g_0& \\
f_2 & f_1 & f_0 &     &        &      &     & g_2 & g_1 & \fbox{$g_0$} \\
\Cline{1-9} 
f_3 & f_2 & f_1 & f_0 &        &      &     & g_3 & g_2 & g_1 \\  
    & f_3 & f_2 & f_1 & \ddots &      &     & \vdots & \vdots & \vdots \\ 
    &     & f_3 & f_2 & \ddots & f_0  &     & \vdots & \vdots & \vdots \\ 
    &     &     & f_3 & \ddots & f_1 & f_0  & \vdots & \vdots & \vdots \\ 
    &     &     &     & \ddots & f_2 & f_1  &   g_n  & g_{n-1}&  g_{n-2} \\ 
   &     &     &     &        & f_3 & f_2  &        & \fbox{$g_n$} & g_{n-1} \\ 
\Cline{1-8} \Cline{10-10}
    &     &     &     &        &     & f_3  &        &        &  g_n   \\
\end{array} \right] \]

First notice that this must have a factor of $f_3$ from the last row.
We see that there are two possibilities for the first column.
Either it is $f_1$ or $f_3$.
If it is $f_1$, then the remainder of the expression is given by 
    $F(m,k-1, k', m'-1)$.
If it is $f_3$, then we see that the second column must contain $f_0$.
After this, the remainder of the expression is given by $-F(m-1, k, k', m'-2)$.
Thus we see that this case will contribute 
    \[ -1\times(F(m, k-1, k', m'-1) - F(m-1, k, k', m'-2)). \]

Here the $-1$ in front comes from the sign of the matrix of the $g_0^2 g_n$.

\item
\[\left[
\begin{array}{cccccccccc}
f_0 &     &     &     &        &      &     & g_0 &    & \\ 
f_1 & f_0 &     &     &        &      &     & g_1 & \fbox{$g_0$}& \\
\Cline{1-8} \Cline{10-10}
f_2 & f_1 & f_0 &     &        &      &     & g_2 & g_1 & \fbox{$g_0$} \\
\Cline{1-9} 
f_3 & f_2 & f_1 & f_0 &        &      &     & g_3 & g_2 & g_1 \\  
    & f_3 & f_2 & f_1 & \ddots &      &     & \vdots & \vdots & \vdots \\ 
    &     & f_3 & f_2 & \ddots & f_0  &     & \vdots & \vdots & \vdots \\ 
    &     &     & f_3 & \ddots & f_1 & f_0  & \vdots & \vdots & \vdots \\ 
  &     &     &     & \ddots & f_2 & f_1  & \fbox{$g_n$} & g_{n-1}&  g_{n-2} \\ 
\Cline{1-7} \Cline{9-10}
    &     &     &     &        & f_3 & f_2  &        &  g_n   &  g_{n-1} \\ 
    &     &     &     &        &     & f_3  &        &        &  g_n   \\
\end{array} \right] \]

With a little work we see that this will contribute 
    \( F(m-1, k, k', m'-2). \)
\end{enumerate}

This combines together to give that 
\begin{eqnarray*} H_0(m,k,k',m') & = & F(m, k, k', m') - F(m, k-1, k', m'-1) \\
                                 &   & + 2 F(m-1, k, k', m'-2).\end{eqnarray*}
By noticing that $F(m,k,k',m') = F(m-1, k, k', m'-2) - F(m,k-1, k', m'-1)
     + F(m,k,k'-1,m')$ we get
\begin{eqnarray*} H_0(m,k,k',m')& = & 2 F(m, k, k', m') + F(m-1, k, k', m'-2) \\
                                &   & - F(m,k,k'-1, m').\end{eqnarray*}
which is the desired result.
\end{proof}

From here we can prove on of the main results which will help us to prove Theorem \ref{thm:mt2}.

\begin{theorem}
\label{thm:tight}
Let $\beta \approx 8.13488$ be the real root of $x^3-18 x^2+110 x-242$, and
    $\alpha \approx 1.83928$ be the real root of $x^3-x^2-x-1$.
Then
    \[H_0(n) = \frac{\beta}{n \pi}  \alpha^n - 
    \OO\left(\frac{\alpha^n}{n^2}\right).\]
\end{theorem}

In order to prove Theorem \ref{thm:tight}, we will find an asymptotic for 
    $H_0(n)$ by maximizing $H_0(m,k,k', m')$ over the real numbers, and then 
    accounting for the error introduced.

\begin{proof}[Proof of Theorem \ref{thm:tight}]
Let us find where $|H_0(m,k,k',m')|$ is maximized.
(Notice that $m'$ is completely determined by $k$ and $m$, and further
    that $n = 3 m + 2 k + k'$).
By writing the factorials as $\Gamma$ functions, and ignoring the $(-1)^k$
    we are maximizing
\[ \hat H(m, k, k') = (3 m + 2k + k')\frac{\Gamma(m+k+k')} {\Gamma(k+1) 
    \Gamma(m+1) \Gamma(k'+1)} \]
subject to the condition 
\[ G(m, k, k') = 3 m + 2 k + k' = n. \]

Thus, to solve for the maximums, we use Lagrange multipliers
   to solve the equations: 
\[ \nabla \hat H = \lambda \nabla G\ \mathrm{and}\ G(m,k,k') = n. \]
Recall that $\Psi(x)$ denotes the digamma function of $x,$ i.e. $\Psi(x)=\frac{\Gamma'(x)}{\Gamma(x)}.$
The latter gives rise to the equations:
\begin{eqnarray*}
3 \lambda & = & 
(3m+2k+k')\frac{\Gamma(m+k+k')}{\Gamma(k+1)\Gamma(m+1)\Gamma(k'+1)}\Psi(k'+k+m)-
\\ &&
  (3m+2k+k')\frac{\Gamma(m+k+k')}{\Gamma(k+1)\Gamma(m+1)\Gamma(k'+1)}\Psi(m+1)+
\\ &&
  3\frac{\Gamma(m+k+k')}{\Gamma(k+1)\Gamma(m+1) \Gamma(k'+1)}\\
2 \lambda & = & 
(3m+2k+k')\frac{\Gamma(m+k+k')}{\Gamma(k+1)\Gamma(m+1)\Gamma(k'+1)}\Psi(k'+k+m)-
\\ &&
  (3m+2k+k')\frac{\Gamma(m+k+k')}{\Gamma(k+1)\Gamma(m+1)\Gamma(k'+1)}\Psi(k+1)+
\\ &&
  2 \frac{\Gamma(m+k+k')}{\Gamma(k+1)\Gamma(m+1) \Gamma(k'+1)}\\
\lambda & = & 
(3m+2k+k')\frac{\Gamma(m+k+k')}{\Gamma(k+1)\Gamma(m+1)\Gamma(k'+1)}\Psi(k'+k+m)-
\\ &&
  (3m+2k+k')\frac{\Gamma(m+k+k')}{\Gamma(k+1)\Gamma(m+1)\Gamma(k'+1)}\Psi(k'+1)+
\\ &&
  \frac{\Gamma(m+k+k')}{\Gamma(k+1)\Gamma(m+1) \Gamma(k'+1)}\\
n  & = & 3 m + 2 k + k'
\end{eqnarray*}

Upon some simplification this becomes:
\begin{eqnarray*}
3 \lambda & = & F(m,k,k') (\Psi(k'+k+m)-\Psi(m+1)+3 /n)  \\
2 \lambda & = & F(m,k,k') (\Psi(k'+k+m)-\Psi(k+1)+2 /n)  \\
\lambda & = & F(m,k,k') (\Psi(k'+k+m)-\Psi(k'+1)+1 /n)  \\
n & = & 3 m + 2 k + k'.
\end{eqnarray*}
By redefining $\lambda$, we get
\begin{eqnarray*}
3 \lambda & = & \Psi(k'+k+m)-\Psi(m+1)+3/n \\
2 \lambda & = & \Psi(k'+k+m)-\Psi(k+1)+2/n \\
\lambda & = & \Psi(k'+k+m)-\Psi(k'+1)+1/n \\
n & = & 3 m + 2 k + k'.
\end{eqnarray*}

If we solve for $\lambda-1/n$ in these equations, and equate them, we
    get the following three equations:
\begin{eqnarray*}
\Psi(k'+k+m)-\Psi(k'+1) & = & \frac{\Psi(k'+k+m)-\Psi(m+1)}{3} \\
\frac{\Psi(k'+k+m)-\Psi(k+1)}{2} & = & \Psi(k'+k+m)-\Psi(k'+1) \\
n & = & 3 m + 2 k + k'.
\end{eqnarray*}

By noticing that $\Psi(n) = \ln(n) + \OO(1/n)$, we can rewrite this as
\begin{eqnarray}
\frac{2}{3} \ln(k' + k + m) - \ln(k' +1 ) + \frac{1}{3}\ln(m+1)   
    & = & \OO\left(\frac{1}{n}\right) \label{eq:ln1} \\
\frac{1}{2}\ln(k' + k + m)   - \ln(k + 1)  + \frac{1}{2}\ln(k' + 1) 
    & = & \OO\left(\frac{1}{n}\right) \label{eq:ln2} \\  
n & = & 3 m + 2 k + k'. \label{eq:n} 
\end{eqnarray}
Here we use the fact that $\OO(k) = \OO(m) = \OO(k') = \OO(n)$.

Now, the question is, what sort of error do we get in the solution of the 
    equations.
For large $k'$, $k$ and $m$, the right hand side is approximately 0, so we can
    find the solution for 0, and then figure out how far off we are.  
Thus we need to find a bound for how quickly the left hand side 
    can change (i.e. derivative), and then figure out how skewed 
    the solution is.

The gradients of the left hand sides are
\[
  \left[ \frac{2}{3 (k' + k + m)}, 
      \frac{2}{3 (k' + k + m)} - \frac{1}{k' + 1},
      \frac{2}{3 (k' + k + m)} + \frac{1}{3 (m + 1)} \right] \]
\[ \left[ \frac{1}{2 ( k' + k + m)} - \frac{1}{2 (k + 1)}, 
      \frac{1}{2 ( k' + k + m)} + \frac{1}{2 (k' + 1)}, 
      \frac{1}{2 ( k' + k + m)} \right]. \]

So we notice that the maximal directional derivatives are $\OO(1/n)$.
This means that the maximal deviation from the actual solution
    is $\OO(1)$.

By solving equations (\ref{eq:ln1}), (\ref{eq:ln2}) and (\ref{eq:n}), where the 
    right hand size is 0 (via Maple \cite{Maple}) 
    and accounting for the $\OO(1)$ term, we can write
\begin{eqnarray*}
m & = & \hat m n + \Delta m \\
k & = & \hat k n + \Delta k \\
k' & = & \hat k' n + \Delta k' \\
\end{eqnarray*}
where $\Delta m$, $\Delta k$ and $\Delta k'$ are all $\OO(1)$, and such
    that $m, k$ and $k'$ are integers, and further that 
\begin{eqnarray*}
\hat m & = & -{\frac {1}{66}}\sqrt [3]{1331+231\sqrt {33}}-1/3{\frac {1}{
               \sqrt [3]{1331+231\sqrt {33}}}}+1/3 \\
\hat k & = & {\frac {1}{66}}\sqrt [3]{3267+627\sqrt {33}}-2{\frac {1}{\sqrt[3]
              {3267+627\sqrt {33}}}} \\
\hat k' & = & {\frac {1}{66}}\sqrt [3]{3267+561\sqrt {33}}+{\frac {1}{\sqrt[3]
              {3267+561\sqrt {33}}}}.
\end{eqnarray*}

We notice that, asymptotically:
\begin{eqnarray*}
& & \hat H(\hat m n + \Delta m,\hat k n + \Delta k,\hat k' n+ \Delta k')  \\
   & = & n \frac{\Gamma((\hat m+\hat k+\hat k')n+\Delta m+\Delta k+\Delta k')} 
      { \Gamma(\hat m n+1 + \Delta m) \Gamma(\hat k n+1 + \Delta k) 
       \Gamma(\hat k' n+1 + \Delta k')}  \\
   & \approx & n \frac{((\hat m+\hat k+\hat k')n)^{\Delta m+\Delta k+\Delta k'}
      \Gamma((\hat m+\hat k+\hat k')n)} 
      { (\hat m n+1)^{\Delta m} \Gamma(\hat m n+1) 
       (\hat k n+1)^{\Delta k} \Gamma(\hat k n+1) 
       (\hat k' n+1)^{\Delta k'} \Gamma(\hat k' n+1) } \\
   & \approx & \frac{(\hat m+\hat k+\hat k')^{\Delta m+\Delta k+\Delta k'}}
      { \hat m^{\Delta m} \hat k^{\Delta k} \hat k'^{\Delta k'}}
      \times n \frac{\Gamma((\hat m+\hat k+\hat k')n)} 
      { \Gamma(\hat m n+1) \Gamma(\hat k n+1) \Gamma(\hat k' n+1) } \\
   & = & \OO(1)n \frac{\Gamma((\hat m+\hat k+\hat k')n)} 
      {\Gamma(\hat k n+1) \Gamma(\hat m n+1) \Gamma(\hat k' n+1) } \\
   & = & \OO(1) \left( 
         \frac{\beta}{\pi n} \alpha^n - \OO\left(\frac{\alpha^n}
         {n^2}\right)\right).
\end{eqnarray*}

Let us consider this $\OO(1)$ term more precisely.
Notice that, using the property that $3 \Delta m + 2 \Delta k + \Delta k' = 0$,
    we have:
\begin{eqnarray*}
 & & \frac{(\hat m+\hat k+\hat k')^{\Delta m+\Delta k+\Delta k'}}
      { \hat m^{\Delta m} \hat k^{\Delta k} \hat k'^{\Delta k'}} \\
& = & \frac{(\hat m+\hat k+\hat k')^{\Delta m+\Delta k-3 \Delta m - 2 \Delta k}}
      {\hat m^{\Delta m} \hat k^{\Delta k} \hat k'^{-3 \Delta m - 2\Delta k}} \\
& = & \frac{(\hat m+\hat k+\hat k')^{- 2\Delta m-\Delta k}}
      { \hat m^{\Delta m} \hat k^{\Delta k} \hat k'^{-3 \Delta m- 2\Delta k}} \\
& = & \frac{(\hat m+\hat k+\hat k')^{- 2\Delta m}
            (\hat m+\hat k+\hat k')^{- \Delta k}}
      { \hat m^{\Delta m} \hat k^{\Delta k} \hat k'^{-3 \Delta m}
        \hat k'^{- 2\Delta k}} \\
& = & \frac{ \hat k'^{3 \Delta m} \hat k'^{ 2\Delta k} }
    { \hat m^{\Delta m} (\hat m+\hat k+\hat k')^{2\Delta m} \hat k^{\Delta k} 
    (\hat m+\hat k+\hat k')^{\Delta k} } \\
& = & \frac{ \hat k'^{3 \Delta m} }{ \hat m^{\Delta m} 
             (\hat m+\hat k+\hat k')^{2\Delta m} }
      \frac{ \hat k'^{ 2\Delta k} } { \hat k^{\Delta k} 
             (\hat m+\hat k+\hat k')^{\Delta k} } \\
& = & \left(\frac{ \hat k'^{3 } }{ \hat m 
             (\hat m+\hat k+\hat k')^{2}}\right)^{\Delta m}
      \left(\frac{ \hat k'^{ 2} } { \hat k
             (\hat m+\hat k+\hat k') }\right)^{\Delta k} \\
& = & 1^{\Delta m} 1^{\Delta k} \\
& = & 1
\end{eqnarray*}
where this last simplification was done via Maple.

So this becomes 
    \[H_0(n) = \frac{\beta}{n \pi}  \alpha^n - \OO\left(\frac{\alpha^n}
     {n^2}\right)\]
    where $\beta$ is the real root of $x^3-18 x^2+110 x-242$, and
    $\alpha$ is the real root of $x^3-x^2-x-1$.
\end{proof}

Theorem \ref{thm:mt2} follows directly from Theorem \ref{thm:tight} 
    and the following Lemma

\begin{lemma}
\label{lem:alpha}
For $n$ sufficiently large, $H_l(n) \leq H_0(n)$.
\end{lemma} 

\begin{proof}
From the comments following the statement of Theorem \ref{thm:main} 
    we see that 
\begin{eqnarray*}
H_l(m, k, k', m') & = & H_l(m, k, k'-1, m') - H_l(m, k-1, k', m'-1)  \\ && 
                       + H_l(m-1, k, k', m'-2).
\end{eqnarray*}

From this it follows that 
\[H_l(n) \leq H_l(n-1) + H_l(n-2) + H_l(n-3) \]

Notice that 
\begin{equation}\label{rec}
H_l(n) = H_{n-l}(n)
\end{equation} 
by considering the resultant with the reciprocal polynomial, 
    namely that ${\rm Res}(f,g)=\pm{\rm Res}(x^3f(1/x),x^ng(1/x))$.
So, we can suppose w.l.o.g. that $l\geq\frac{n}{2}$.
We write this as
\begin{eqnarray*}
H_l(n) & \leq & 1 \times H_l(n-1) + 1 \times H_l(n-2) + 1 \times H_l(n-3)\\
       & :=   & A_1 H_l(n-1) + B_1 H_l(n-2) + C_1 H_l(n-3)\\
       & \leq & (A_1 + B_1) H_l(n-2) + (A_1 + C_1) H_l(n-3) + A_1 H_l(n-4)\\
       & :=   & A_2 H_l(n-2) + B_2 H_l(n-3) + C_2 H_l(n-4) \\
       & \vdots & \\
       & \leq & A_{n-l-2} H_l(l+2) + B_{n-l-2} H_l(l+1) + C_{n-l-2} 
                H_l(l) \\
       & = & A_{n-l-2} H_2(l+2) + B_{n-l-2} H_{1}(l+1) + C_{n-l-2} 
                H_{0}(l), \\
\end{eqnarray*}
where the last equality holds because of (\ref{rec}).
The numbers $A_m, B_m$ and $C_m$ satisfy linear recurrence relationships.
Namely, we have that $A_m = A_{m-1} + B_{m-1}, B_m = A_{m-1} + C_{m-1}$ and 
    $C_m = A_{m-1}$.
This simplifies to $A_1 = 1, A_2 = 2, A_3 = 4, A_m = A_{m-1} + A_{m-2} 
    + A_{m-3}$, and further that $B_m = A_{m-1} + A_{m-2}$ and $C_m = A_{m-1}$.

Solving this gives
$A_m = c \alpha^m + c_1 \alpha_1^m + c_2 \alpha_2^m$, where
$\alpha$ is the real root of $x^3-x^2-x-1$, and $\alpha_i$ are its 
    conjugates.
Further $c$ is the real root of $44 x^3 - 44 x^2 + 12 x -1$ and
     $c_1$ and $c_2$ are its conjugates.

Numerically 
\begin{eqnarray*}
c & \approx & .6184199224 \\
c_1 & \approx & .1907900391 + .01870058339 i \\
c_2 & \approx & .1907900391 - .01870058339 i
\end{eqnarray*}

For $m \geq 3$, this gives us by the triangle inequality,
    $A_m \leq 0.7 \alpha^m$.
Similarly, for $m \geq 5$ we get that 
\[B_m = A_{m-1} + A_{m-2} \leq \alpha^m (0.7 / \alpha + 0.7/\alpha^2)
                         \leq 0.6 \alpha^m \]
    and for $m \geq 4$ we get that
\[C_m = A_{m-1} \leq \alpha^m (0.7/\alpha) \leq 0.4 \alpha^m \]

Now, we have already shown that 
\[H_0(n) = \frac{\beta}{\pi n}\alpha^n-\OO\left(\frac{\alpha^n}{n^2}\right)\]
    where $\beta = 8.13488$ (Theorem \ref{thm:tight}).

Using the same method, we can show that
\[H_l(n) = \frac{\beta_l}{\pi n}\alpha^n-\OO\left(\frac{\alpha^n}{n^2}\right)\]
for $l$ from $0$ to $6$ where
\begin{eqnarray*}
\beta_0 & = & 8.13488 \\
\beta_1 & = & 3.71205 \\
\beta_2 & = & 0.92093 \\
\beta_3 & = & 1.01680 \\
\beta_4 & = & 0.31597 \\
\beta_5 & = & 0.01923 \\
\beta_6 & = & 0.05956 
\end{eqnarray*}
So, $H_l(n)\leq H_0(n)$ if $n-6\leq l\leq n$ (this is again due to (\ref{rec})).
Suppose now that $l\leq n-7.$ Then $n-l-2\geq5$ and all the bounds computed above for $A_m,B_m,C_m$ hold.
So, we have, for large $n$, 
\begin{equation}
\label{huno} \begin{array}{rcl}
H_l(n) 
& \leq & A_{n-l-2} H_2(l+2) + B_{n-l-2} H_{1}(l+1) + C_{n-l-2} H_{0}(l)  \\
& \leq & 0.7 \alpha^{n-l-2} \left(\frac{\beta_2}{\pi(l+2)}\alpha^{l+2}
        -\OO(\frac{\alpha^{l+2}}{(l+2)^2})\right) \\
&      & + 0.6 \alpha^{n-l-2} \left(\frac{\beta_1}{\pi(l+1)}\alpha^{l+1}
        -\OO(\frac{\alpha^{l+1}}{(l+1)^2})\right) \\
&      & +0.4 \alpha^{n-l-2} \left(\frac{\beta_0}{\pi l}\alpha^{l}
        -\OO(\frac{\alpha^{l}}{(l)^2})\right) \\
& \leq &  0.7 \alpha^{n-l-2} \frac{\beta_2}{\pi(l+2)}\alpha^{l+2}
        + 0.6 \alpha^{n-l-2} \frac{\beta_1}{\pi(l+1)}\alpha^{l+1} \\
&      & + 0.4 \alpha^{n-l-2} \frac{\beta_0}{\pi l}\alpha^{l} \\
&    = &  0.7 \frac{\beta_2}{\pi(l+2)}\alpha^{n}
        + 0.6 \frac{\beta_1}{\pi(l+1)}\alpha^{n-1}
        +0.4  \frac{\beta_0}{\pi l}\alpha^{n-2}.
\end{array}
\end{equation}
The last expression of (\ref{huno}) is maximal when $l$ is minimal, i.e. 
    $l=n/2.$ 
So, for large $n$, we get that $H_l(n)$ is bounded above by
\begin{eqnarray*}
H_l(n) & \leq & 
            0.7 \frac{\beta_2}{\pi(n/2+2)}\alpha^{n}
           +0.6 \frac{\beta_1}{\pi(n/2+1)}\alpha^{n-1}
           +0.4 \frac{\beta_0}{\pi n/2}\alpha^{n-2} \\
   & \leq & 0.7 \frac{\beta_2}{\pi(n/2)}\alpha^{n}
           +0.6 \frac{\beta_1}{\pi(n/2)}\alpha^{n-1}
           +0.4 \frac{\beta_0}{\pi n/2}\alpha^{n-2} \\
   & \leq & 2 \left(0.7 \times \beta_2
           +0.6 \frac{\beta_1}{\alpha}
           +0.4 \frac{\beta_0}{\alpha^2}\right)\frac{\alpha^n}{\pi n}\\
   & = & \frac{5.6348}{\pi n} \alpha^n
\end{eqnarray*}
This expression is bounded above by $H_0(n)=\frac{\beta_0}{\pi n}\alpha^n
    -\OO(\frac{\alpha^n}{n^2})$ for large values of $n,$ 
    which gives the desired result.
\end{proof}

\bigskip
Now we are ready for the proof of our main result.
\begin{proof}[Proof of Theorem \ref{thm:mt2}]
Due to Theorem \ref{thm:tight}, we will be done if we show that, for $n\gg0,\,H(n)=H_0(n).$
As it was showed in Lemma \ref{lem:alpha}, it turns out that $H_0(n)=\max_{0\leq l\leq n} H_l(n)$ if $n\gg0.$ 
As explained at the beginning of this section, 
notice that if $H(n) > H(n-1)$ then $H(n) = \max_l H_l(n),$ so we only have to prove that for infinite values of $N,$ we have 
$H(N)>H(N-1).$
\par Suppose this is not the case, then $H(N)$ is bounded as $N\to\infty,$ and this is a contradiction with Theorem \ref{thm:tight} which
says that $H(N)\geq H_0(N)_{N\to\infty}\to+\infty.$
\par
So pick $N$ such that $H(N) > H(N-1)$, and sufficiently large such that 
    $H(N) = H_0(N) \geq \max_l H_l(N)$ (Lemma \ref{lem:alpha}) and 
    $H(N+1) \geq H_0(N+1) > H_0(N)$.
Hence by induction for all $m \geq N$ we have that $H(m) > H(m-1)$ 
    and $H(m) = H_0(m)$.

\end{proof}

It should be pointed out that experimentally, $H(n) > H(n-1)$ for all $n$
    and $H(n) = H_0(n)$ for all $n \geq 18$.

\section{Conclusions and comments}
\label{sec:conc}
In this paper we give a precise description for $H(\R(f,g))$ where
    $f$ is a quadratic polynomial, and tight asymptotics 
    when $f$ is a cubic polynomial.
The methods used in this paper should be extendible to the case of $f$ being
    a polynomial of fixed degree $m$.
In particular, most of Section \ref{sec:tri} is done constructively,
    and can be extended to arbitrary $m$.
So we can most likely find bounds such as $H(n) \leq \OO(\alpha^n)$
    for arbitrary fixed $m$, and $\alpha$ dependent on $m$.
It would be interesting and worthwhile to do this.

Let $g(x) = g_0 + \cdots + g_n x^n$ be a degree $n$ polynomial.
As a result of Lemma \ref{lem:alpha} we proved that for 
    sufficiently large $n$ that 
\[ H(\R(f_0 + \cdots + f_3 x^3,g))
  = H\left(\R(f_0+\cdots+f_3 x^3, g_0 + g_n x^n)\right).\]
(Experimentally, this appears to be true for $n \geq 18$.)
Notice that if $\deg(f)=2,$ for $n \geq 3$:
\[
H(\R(f_0 + f_1 x + f_2 x^2, g)) = H(\R(f_0 + f_1 x + f_2 x^2,g_0+g_n x^n)).
\]
It is trivial to see that in the linear case:
\begin{eqnarray*}
H(\R(f_0 + f_1 x, g)) & = & H(\R(f_0 + f_1 x,g_0+g_n x^n))  \\
                      & ( = & 1).
\end{eqnarray*}

It is reasonable to conjecture that
\begin{conjecture}
For fixed $m$, and $g(x) = g_0 + \cdots + g_n x^n$, for sufficiently 
    large $n$ (dependent on $m$), 
\[ H(\R(f_0 + \dots + f_m x^m, g)) = 
   H(\R(f_0 + \dots + f_m x^m, g_0 + g_n x^n)). \]
\end{conjecture}

There is some computational evidence to support this conjecture.

\section{Acknowledgments}

We are grateful to Mart\'{\i}n Sombra for providing us updated references concerning the state of the art of the computation of heights of resultants.
We are also grateful to Teresa Krick for helpful comments on a preliminary version of this paper.

\section{Appendix: Tables}
\label{sec:append}

\begin{table}[H]
\begin{tabular}{|cl|cl|cl|}
\hline
$A_n$ & $n$ & $A_n$ & $n$ & $A_n$ & $n$ \\
\hline
1& 3,4 & 10& 34,35,36,37 & 19& 67,68,69,70\\
2&5,6,7,8 & 11& 38,39,40,41 & 20&71,72,73\\
3& 9,10,11,12& 12& 42,43,44 & 21&74,75,76,77\\
4& 13,14,15 & 13& 45,46,47,48 & 22&78,79,80,81\\
5& 16,17,18,19 & 14& 49,50,51,52& 23&82,83,84\\
6& 20,21,22,23& 15& 53,54,55& 24&85,86,87,88\\
7& 24,25,26 & 16& 56,57,58,59& 25&89,90,91\\
8& 27,28,29,30 & 17& 60,61,62& 26&92,93,94,95\\
9& 31,32,33& 18& 63,64,65,66& 27&96,97,98,99 \\
\hline
\end{tabular}
\caption{Values of $A_n$ (Theorem \ref{thm:mt1}, page \pageref{thm:mt1})}
\label{tab:An}
\end{table}

\begin{table}[H]
\begin{tabular}{|cl|cl|cl|}
\hline
   n   &  Maximum at $H_l$ &
   n   &  Maximum at $H_l$ &
   n   &  Maximum at $H_l$ \\
\hline
   1   &    $H_0$   & 8    &    $H_0$   & 15   &    $H_3$   \\
   2   &    $H_1$   & 9    &    $H_3$   & 16   &    $H_3$   \\
   3   &    $H_0$   & 10   &    $H_3$   & 17   &    $H_3$   \\

   4   &    $H_1$   & 11   &    $H_0$   & 18   &    $H_0$   \\
   5   &    $H_1$ and $H_2$ 
                    & 12   &    $H_0$   & 19   &    $H_0$   \\
   6   &    $H_3$   & 13   &    $H_3$   & \vdots&   \vdots  \\ 
   7   &    $H_3$   & 14   &    $H_3$   & 72   &    $H_0$   \\
\hline
\end{tabular}
\caption{Maximal $H_l$ value}
\label{tab:Hl}
\end{table}

\begin{table}[H]
\[
\begin{array}{|rcp{3.5 in}|}
\hline
H_0(m,k,k',m') & = & $F(m-1,k,k',m'-2)-F(m,k,k'-1,m')+2F(m,k,k',m')$\\ 

H_1(m,k,k',m') & = & $2  F(m-1,k,k'-1,m'-1)- F(m,k-1,k'-1,m')+2  F(m,k-1,k',m')
    -3  F(m-1,k,k',m'-1)$\\ 

H_2(m,k,k',m') & = & $2  F(m-1,k,k'-2,m')-4  F(m-1,k,k'-1,m') 
    -F(m-2,k-1,k',m'-3)
    -3  F(m-2,k,k',m'-2)+ F(m-1,k-2,k',m'-2)- F(m,k-2,k'-1,m')
    +2  F(m,k-2,k',m')$\\ 

H_3(m,k,k',m') & = & $-2  F(m-2,k,k'-2,m'-1)+3  F(m-1,k-1,k'-2,m')
    -6  F(m-1,k-1,k'-1,m')
    + F(m-3,k,k',m'-3)+5  F(m-2,k,k',m'-1)-2  F(m-2,k-1,k',m'-2)
    - F(m-2,k-2,k',m'-3)+ F(m-1,k-3,k',m'-2)- F(m,k-3,k'-1,m')
    +2  F(m,k-3,k',m')$\\

H_4(m,k,k',m') & = & $-2  F(m-5,k,k',m'-6)- F(m-4,k,k',m'-4)
    +3  F(m-3,k-1,k'-1,m'-3) 
    -9  F(m-2,k-2,k'-1,m'-2)+ F(m-2,k-3,k',m'-3)-7  F(m-2,k-2,k',m'-2)
    +13  F(m-3,k-1,k',m'-3)+6  F(m-3,k,k'-2,m'-2)+2  F(m-2,k,k'-3,m')
    + F(m-1,k-4,k',m'-2)- F(m,k-4,k'-1,m')+2  F(m,k-4,k',m')
    +4  F(m-1,k-2,k'-2,m')-8  F(m-1,k-2,k'-1,m')$\\

H_5(m,k,k',m') & = & $2  F(m-3,k,k'-3,m'-1)+18  F(m-3,k-1,k'-2,m'-2)
    -7  F(m-3,k,k'-2,m'-1)
    +12  F(m-4,k-1,k'-1,m'-4)-13  F(m-4,k,k'-1,m'-3)- F(m-5,k-1,k',m'-6)
    -3  F(m-5,k,k',m'-5)+5  F(m-2,k-1,k'-2,m')+2  F(m-1,k-5,k',m'-2)
    + F(m,k-5,k',m')- F(m,k-6,k',m'-1)+5  F(m-1,k-4,k'-1,m'-1)
    -5  F(m-1,k-3,k'-1,m')-15  F(m-2,k-4,k',m'-3)-25  F(m-2,k-3,k',m'-2) 
    +10  F(m-3,k-2,k'-1,m'-3)+15  F(m-4,k-2,k',m'-5)$ \\
\hline
\end{array}
\]
\caption{A table of $H_l(m,k,k',m')$ values, (Theorem \ref{thm:main}, 
         page \pageref{thm:main})}
\label{tab:recurrence}
\end{table}

\end{document}